\documentclass{elsart}

\makeatletter
\long\def\@maketablecaption#1#2{\@tablecaptionsize
    \global \@minipagefalse
    \hbox to \hsize{\parbox[t]{\hsize}{\centering #1 \\ #2}}}

\makeatother

\usepackage{graphics}
\usepackage{graphicx}      
\usepackage{natbib}        
\usepackage{epsfig} 
\usepackage{psfrag} 
\usepackage{amsmath}
\usepackage{amssymb}
\usepackage{mathrsfs}
\usepackage{eufrak}

\begin{document}
\begin{frontmatter}

\title{Stable Hierarchical Model Predictive Control Using an Inner Loop Reference Model and $\lambda$-Contractive
Terminal Constraint Sets}

\author[First]{Chris Vermillion}
\author[Second]{Amor Menezes}
\author[Third]{Ilya Kolmanovsky}

\address[First]{Altaeros Energies, Boston, MA 02110 (e-mail: chris.vermillion@altaerosenergies.com)}
\address[Second]{California Institute for Quantitative Biosciences, University of California, Berkeley, Berkeley, CA 94704 (e-mail: amenezes@berkeley.edu)}
\address[Third]{Department of Aerospace Engineering, University of Michigan, Ann Arbor, MI 48109 (e-mail: ilya@umich.edu)}

\begin{abstract}                
This paper presents a hierarchical model predictive control (MPC)
framework that is presented in \cite{hierarchical_vermillion} and \cite{hierarchical_vermillion_journal}, along with
proofs that were omitted in the aforementioned works.  The
method described in this paper differs
significantly from previous approaches to guaranteeing overall stability, which
have relied upon a multi-rate framework where the inner loop (low level) is updated
at a faster rate than the outer loop (high level), and the inner loop must reach
a steady-state within each outer loop time step.  In contrast,
the method proposed in this paper is aimed at stabilizing the origin of an error system characterized
by the difference between the inner loop state and the state specified by a full-order
reference model.  This makes the method applicable to systems with reduced levels
of time scale separation.  This paper reviews the fundamental results of \cite{hierarchical_vermillion} and \cite{hierarchical_vermillion_journal} and presents proofs that were omitted due to space limitations.
\end{abstract}


\end{frontmatter}

\section{Introduction}\label{intro}
\begin{figure}[thpb]
\centering
\includegraphics[width=3.0in]{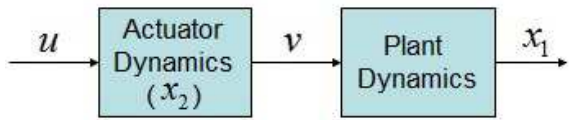}
\caption{Block diagram of the actuator/plant cascade considered in this work.} \label{hierarchical_ol}
\end{figure}

This paper focuses on a two-layer inner loop/outer loop hierarchical control
structure where the ultimate objective
is to stabilize the overall system.  The actuator and plant represent a cascade, depicted in Fig. \ref{hierarchical_ol}, wherein
an actuator output, denoted by $v$, characterizes
an overall force, moment, or generalized effect produced by the actuators, and is
referred to as a \emph{virtual control input}. In the hierarchical control
strategy, an outer loop controller sets a desired
value for this virtual control input, denoted by $v_{des}$, and it is the responsibility of the inner loop to generate control
inputs $u$ that drive $v$ close to $v_{des}$.

This control approach is employed in a number of automotive, aerospace, and marine applications,
such as \cite{ca2}, \cite{ca3}, \cite{ca4}, \cite{ca5}, and \cite{ca6}. Hierarchical
control has become commonplace in industrial applications, as it offers two key advantages
over its centralized counterpart:
\begin{enumerate}
  \item Plug-and-play integration of new design features (for example, a new inner loop) without
  requiring a complete system redesign;
  \item Reduction in overall computational complexity, in terms of the number of inputs and/or states
  considered by each controller.
\end{enumerate}
The use of MPC for constrained hierarchical control has been a natural choice in instances when constraint satisfaction was critical and/or multiple control objectives were traded off. Only recently, however, has an effort been made to provide theoretical stability guarantees for the hierarchical system. In many recent papers, including \cite{hierarchical_app1}, \cite{hierarchical1}, \cite{hierarchical2}, \cite{hierarchical_review}, and \cite{hierarchical3}, the inner
loop is updated at a faster rate than the outer loop, and the inner loop is designed to reach a steady-state, wherein $v=v_{des}$,
within a single outer loop time step.  This strategy represents an
effective way of guaranteeing stability under large time-scale separation, but numerous systems, including those described in \cite{ca2}, \cite{ca3}, and \cite{ca4} (which address a flight control application), and \cite{ca6}, \cite{ca7}, and \cite{ca8} (which address a thermal management system), do not exhibit such a demonstrable time scale separation.

Our approach differs from that of \cite{hierarchical1}, \cite{hierarchical2}, and
\cite{hierarchical3} in that it drives the inner loop states to those of a \emph{reference model} rather than to the steady state values corresponding to $v_{des}$. Our stability formulation relies on
$\lambda$-contractive terminal constraint sets for the outer and inner loop, in addition to rate-like constraints
that ensure that the optimized MPC trajectories do not vary too much from one instant
to the next.  The contractive nature of the terminal constraint allows MPC-optimized control input trajectories to vary from one time step and the next.

The work presented in this paper is an extension of an original IFAC conference
paper (\cite{hierarchical_vermillion}) and builds upon it through two key mechanisms:
\begin{itemize}
  \item Allowance for inexact (approximate) inner loop reference model matching;
  \item Greater flexibility in the decay of contraction rates within the MPC optimization.
\end{itemize}

\section{Problem Statement}\label{problem}
In this paper, we consider two interconnected systems, as depicted in Fig. \ref{hierarchical_basic}, whose dynamics in discrete
time are given by:
\begin{eqnarray}\label{ssr_d}
\nonumber x_{1}(k+1)&=&A_{1}x_{1}(k)+B_{1}v(k), \\
x_{2}(k+1)&=&A_{2}x_{2}(k)+B_{2}u(k), \\
\nonumber v(k) &=& Cx_{2}(k),
\end{eqnarray}
\noindent where $v\in\mathbb{R}^{q}$ represents the virtual control
input, $x_{1} \in
\mathbb{R}^{n_{1}}$ represents the plant states, which are driven by
the virtual control input, $v$, whereas $x_{2}\in
\mathbb{R}^{n_{2}}$ represents the actuator states, which are driven
by the real control inputs, $u\in\mathbb{R}^{p}$, where $p \geq q$. The
control inputs, $u$ are subject to a saturation constraint set $U$, such
that $u(k) \in U$ at every time instant.  We assume that:
\begin{itemize}
  \item \emph{Assumption 1:} The pair $(A_{1},B_{1})$ is stabilizable.
  \item \emph{Assumption 2:} The pair $(A_{2},B_{2})$ is controllable.
  \item \emph{Assumption 3:} Without loss of generality, the actuator dynamics of (\ref{ssr_d})
are written in block controllable canonical form (CCF) described in \cite{miso_ccf}.
\end{itemize}
The assumption of stabilizability is clearly essential to any problem
whose objective is stabilization of the origin. The stronger assumption of inner loop controllability (Assumption 2) allows us to generate an
inner loop error system and control law that appropriately places all
of the poles of the closed inner loop to satisfy the reference model specifications.
%
%

\section{Control Design Formulation}\label{control_design}
Our approach relies on the design of an inner loop
\emph{reference model}, which describes the ideal input-output behavior from $v_{des}$ to $v$. We will proceed to derive an error system describing the difference between the inner loop and reference model states, and we will show how closed-form control laws can be used to achieve exact or sufficiently accurate reference model matching near the origin of this system, ultimately resulting in local stability of the overall system. MPC is used to enlarge the region of attraction of the overall system to include states under which the closed-form control laws hit saturation constraints.

The specific control algorithm incorporates both outer and inner loop terminal constraint sets, wherein closed-form control laws are used to achieve reference model matching (or approximate matching) behavior. Farther from the origins of the outer and inner loops, MPC is used to drive the system into these constraint sets, explicitly accounting for saturation constraints.

\subsection{Reference Model Design and Assumptions}
This reference model is given by:
\begin{eqnarray}\label{ref_model}
x_{f}(k+1) &=& A_{f}x_{f}(k)+B_{f}v_{des}(k), \\
\nonumber v_{des}^{f}(k) &=& Cx_{f}(k),
\end{eqnarray}
\noindent where $x_{f}\in
\mathbb{R}^{n_{2}}$, $v_{des} \in \mathbb{R}^{q}$, and $v_{des}^{f} \in \mathbb{R}^{q}$.
We assume that:
\begin{itemize}
  \item \emph{Assumption 4:} The reference model is stable, i.e., $\|\bar{\lambda}_{i}(A_{f}))\| < 1, \forall i$ ($\bar{\lambda}_{i}$ represents the $i^{th}$ eigenvalue of $A_{f}$);
  \item \emph{Assumption 5:} The reference model does not share any zeros with unstable poles
of $A_{1}$.
\end{itemize}
\noindent Because Assumptions 4 and 5 are on the reference model, which is freely chosen by the control system designer, they do not restrict the applicability of the proposed control design.

We use the reference model to analyze the closed-loop behavior of the inner loop through the following error system:
\begin{eqnarray}\label{il_error_dynamics}
\nonumber \tilde{x}(k+1) &=& A_{2}\tilde{x}(k)+(A_{2}-A_{f})x_{f}(k)+B_{2}u(k) \\
&& -B_{f}v_{des}(k), \\
\nonumber \tilde{v}(k) &=& C\tilde{x}(k),
\end{eqnarray}
\noindent where $\tilde{x}(k) = x_{2}(k)-x_{f}(k)$ and it follows that $\tilde{v}(k)=v(k)-v_{des}^{f}(k)$.
For notational convenience
throughout the paper, because the reference model is embedded in the outer loop, we will introduce
the augmented outer loop state, $x_{1}^{aug} \triangleq \left[\begin{array}{cc}x_{1}^{T} & x_{f}^{T}\end{array} \right]^{T}$,
which results in augmented outer loop dynamics given by:
\begin{eqnarray}\label{ol_aug_dynamics}
x_{1}^{aug}(k+1) &=& A_{1}^{aug}x_{1}^{aug}(k) + B_{1}^{aug}\tilde{v}(k) \\
\nonumber && + B_{f}^{aug}v_{des}(k)
\end{eqnarray}
\noindent where:
\begin{eqnarray}
A_{1}^{aug} &=& \left[
\begin{array}{cc} A_{1} & B_{1}C \\
0 & A_{f} \end{array}
\right], \\
\nonumber B_{1}^{aug} &=& \left[\begin{array}{cc} B_{1}^{T} & 0 \end{array}\right]^{T}, \\
\nonumber B_{f}^{aug} &=& \left[\begin{array}{cc} 0 & B_{f}^{T} \end{array}\right]^{T}.
\end{eqnarray}

\subsection{Model Predictive Control Framework}
\begin{figure}[thpb]
\centering
\includegraphics[width=3.5in]{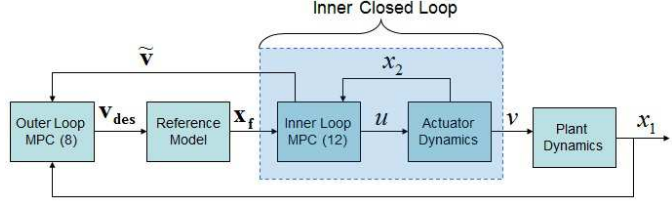}
\caption{Block diagram of the hierarchical control strategy that is implemented when model predictive control is active. In this scenario, N-step
predictions of interconnection variables (shown in \textbf{bold}) are passed between the inner and
outer loop optimizations.} \label{hierarchical_MPC}
\end{figure}

\begin{figure}[thpb]
\centering
\includegraphics[width=3.5in]{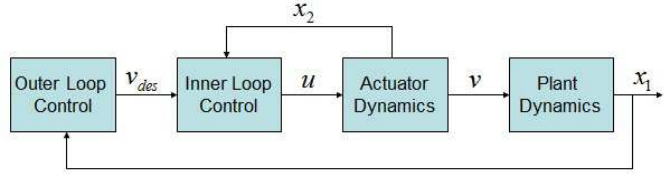}
\caption{Block diagram of the hierarchical control strategy that is implemented under terminal control laws.} \label{hierarchical_basic}
\end{figure}

An MPC optimization is carried out whenever the outer \emph{or} inner loop states are outside of predetermined
$\lambda$-contractive terminal constraint sets $G_{1}$ and $G_{2}$ respectively. A closed-form terminal control law is active once the inner \emph{and} outer loop states have reached the terminal sets. The block diagram of the closed-loop system when MPC is active
is given in Fig. \ref{hierarchical_MPC}, whereas the closed-loop system under closed-form terminal control laws
conforms to the block diagram of Fig. \ref{hierarchical_basic}.

Whenever the MPC optimization is carried out, an optimal control trajectory is computed for an $N$ step prediction
horizon, along with a corresponding state trajectory. The outer loop virtual control and state trajectories are given by:
\begin{eqnarray}\label{sequences}
\nonumber \mathbf{v_{des}}(k) &=& \left[\begin{array}{ccc}\mathbf{v_{des}}(k|k) & \ldots & \mathbf{v_{des}}(k+N-1|k) \end{array} \right], \\
\mathbf{x_{1}^{aug}}(k) &=& \left[\begin{array}{ccc}\mathbf{x_{1}^{aug}}(k|k) & \ldots & \mathbf{x_{1}^{aug}}(k+N|k) \end{array} \right].
\end{eqnarray}
\noindent The inner loop control and state trajectories are given by:
\begin{eqnarray}\label{sequences_inner}
\nonumber \mathbf{u}(k) &=& \left[\begin{array}{ccc}\mathbf{u}(k|k) & \ldots & \mathbf{u}(k+N-1|k) \end{array} \right], \\
\mathbf{\tilde{x}}(k) &=& \left[\begin{array}{ccc}\mathbf{\tilde{x}}(k|k) & \ldots & \mathbf{\tilde{x}}(k+N|k) \end{array} \right].
\end{eqnarray}
\noindent The notation $(i|k)$ denotes the chosen/predicted value of a variable at step $i$ when the optimization is carried
out at time $k$ ($k \leq i$).

The mathematical description of the outer loop control law is:
\begin{eqnarray}\label{mpc_outer_law}
\nonumber v_{des}(k) &=& \left\{ \begin{array}{ccccc} -K_{1}x_{1}^{aug}(k) & \text{if} & x_{1}^{aug}(k) \in G_{1},  \tilde{x}(k) \in G_{2}, \\
 \mathbf{v_{des}^{o}}(k|k) & , & \text{otherwise} \end{array}\right.
\end{eqnarray}
\noindent Here, $K_{1}$ is the terminal control gain and $\mathbf{v}_{des}^{o}(k)$ is the optimized control input sequence from the outer loop MPC optimization, given by:
\begin{equation}\label{mpc_outer_opt}
\mathbf{v_{des}^{o}}(k) = \arg \underset{\mathbf{v_{des}}\in \mathbf{V_{des}}}{\min} J_{1}(\mathbf{v_{des}}(k)|x_{1}^{aug}(k),\mathbf{\tilde{v}}(k-1)),
\end{equation}
subject to the dynamics of (\ref{ol_aug_dynamics}) and constraints:
\begin{eqnarray}\label{mpc_outer_constraints}
\nonumber \mathbf{x_{1}^{aug}}(k+N-1|k) & \in & G_{1}, \\
\mathbf{x_{1}^{aug}}(k+N|k) & \in & \lambda_{1}G_{1}, \\
\nonumber \|\mathbf{v_{des}}(k+i|k)-\mathbf{v_{des}^{o}}(k+i|k-1)\| & \leq & (\delta_{v_{des}}^{max})\beta^{\min(k,N_{1}^{*})}, \\
\nonumber && i=0 \ldots N-2,
\end{eqnarray}
\noindent and cost function:
\begin{eqnarray}\label{mpc_outer_cost}
\nonumber J_{1}(\mathbf{v_{des}}(k)|x_{1}^{aug}(k),\mathbf{\tilde{v}}(k-1)) &=& \sum_{i=k}^{k+N-1}g_{1}(\mathbf{x_{1}^{aug}}(i|k), \\
&& \mathbf{v_{des}}(i|k)).
\end{eqnarray}
\noindent Here, $\lambda_{1}$, $\delta_{v_{des}}^{max}$, $\beta$, and $N_{1}^{*}$ are design
parameters, which are summarized in Table \ref{MPC_parameters}. $\mathbf{V_{des}}$ is the
set of all feasible $\mathbf{v_{des}}$ trajectories. For the results in this paper, there are no restrictions to
the form of the stage cost, $g_{1}(\mathbf{x_{1}^{aug}}(i|k),\mathbf{v_{des}}(i|k))$.
The mathematical description of the inner loop control law is:
\begin{eqnarray}\label{mpc_inner_law}
u(k) &=& \left\{\begin{array}{ccccc} u_{t}(k) & if & x_{1}^{aug}(k) \in G_{1},  \tilde{x}(k) \in G_{2} \\
\nonumber \mathbf{u^{o}}(k|k) & , & otherwise \end{array}\right.
\end{eqnarray}
\noindent where
\begin{equation}\label{inner_terminal}
u_{t}(k) = K_{21}v_{des}(k) - K_{22}x_{2}(k).
\end{equation}
Here, $\mathbf{u^{o}}(k)$ is the optimized control input sequence from the inner loop
MPC optimization, given by:
\begin{equation}\label{mpc_inner_opt}
\mathbf{u^{o}}(k) = \arg \underset{\mathbf{u}(k) \in \mathbf{U}}{\min} J_{2}(\mathbf{u}(k)|\tilde{x}(k),\mathbf{x_{f}}(k)),
\end{equation}
subject to the dynamics of (\ref{il_error_dynamics}) and constraints:
\begin{eqnarray}\label{mpc_inner_constraints}
\nonumber \mathbf{\tilde{x}}(k+N|k) & \in & \lambda_{2}G_{2}, \\
\|\mathbf{u}(k+i|k)-\mathbf{u^{o}}(k+i|k-1)\| & \leq & (\delta_{u}^{max})\beta^{\min(k,N_{2}^{*})}, \\
\nonumber && i=0 \ldots N-2 \\
\nonumber \mathbf{u}(k+i|k) &\in& U, i=0 \ldots N-1,
\end{eqnarray}
\noindent where $U$ reflects the actuator saturation limits of $u$ and $\mathbf{U}$ is the
set of all feasible control input $\mathbf{u}$ trajectories. The inner loop cost function is given by:
\begin{equation}\label{mpc_inner_cost}
J_{2}(\mathbf{u}(k)|\tilde{x}(k),\mathbf{x_{f}}(k)) = \sum_{i=k}^{k+N-1}g_{2}(\mathbf{\tilde{x}}(i|k),\mathbf{u}(i|k)).
\end{equation}
Here, $\lambda_{2}$, $\delta_{u}^{max}$, $\beta$, and $N_{2}^{*}$ are design
parameters, which are summarized in Table \ref{MPC_parameters}. As with the outer loop, there
are no restrictions to the form of the stage cost, $g_{2}(\mathbf{\tilde{x}}(i|k),\mathbf{u}(i|k))$.

The terms $\beta^{\min(k,N_{1}^{*})}$ and $\beta^{\min(k,N_{2}^{*})}$ impose the requirement that
trajectories $\mathbf{v_{des}}(k)$ and $\mathbf{u}(k)$ calculated at any two subsequent time steps must be sufficiently close to each other, and that the required proximity of trajectories decrease over time, until $k=N_{1}^{*}$ and $k=N_{2}^{*}$, respectively.  Formulas for the required values for $N_{1}^{*}$ and $N_{2}^{*}$ are given in the proof of Proposition \ref{convergence}; required values depend on the contraction rates $\lambda_{1}$ and $\lambda_{2}$, system dynamics, horizon length ($N$), and $\beta$.

\begin{figure}[thpb]
\centering
\includegraphics[width=3.5in]{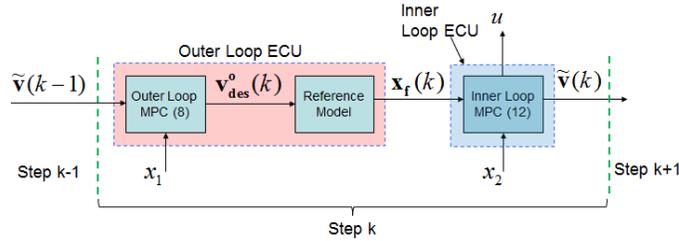}
\caption{Sequence of operations when MPC is active, including the physical locations
where the operations occur.} \label{MPC_sequence}
\end{figure}

Key MPC design parameters, including those that are essential for our stability formulation, are provided in Table \ref{MPC_parameters}. Fig. \ref{MPC_sequence} provides a graphical depiction
of the sequence of operations that occur in a single time instant when MPC is active.

Computationally, the outer loop
MPC must consider $n_{1}+n_{2}$ states and $q$ control inputs, whereas the inner loop must consider
$n_{2}$ states and $p$ control inputs. Both optimizations are individually computationally simpler than their centralized
counterparts, which must consider $n_{1}+n_{2}$ states and $p$ control inputs. The resulting computational simplification
can be especially significant when the algorithm is applied to systems with complex outer loops ($n_{1} >> n_{2}$) and
several actuators for a given virtual control ($p >> q$), which is commonplace in industry.
\begin{table}
  \centering
  \normalsize
  \caption{Key MPC Design Parameters}\label{MPC_parameters}
  \begin{tabular}{|c|c|}
  \hline
  Parameter & Description \\
  \hline
  $G_{1}$ & outer loop terminal constraint set \\
  $G_{2}$ & inner loop terminal constraint set \\
  $K_{1}$ & outer terminal control gain matrix \\
  $K_{21},K_{22}$ & inner terminal control gain matrices \\
  $\lambda_{1}$ & outer contraction rate ($<1$) \\
  $\lambda_{2}$ & inner contraction rate ($<1$) \\
  $\delta_{v_{des}}^{max}$ & rate-like constraint on outer loop MPC \\
  $\delta_{u}^{max}$ & rate-like constraint on inner loop MPC \\
  $\beta$ & any scalar that is $<1$ \\
  $N_{1}^{*}$ & maximum steps until convergence to $G_{1}$ \\
  $N_{2}^{*}$ & maximum steps until convergence to $G_{2}$ \\
  \hline
  \end{tabular}
\end{table}

\section{Deriving Terminal Control Laws and $\lambda$-Contractive Terminal Constraint Sets}\label{lambda_contractive}
In this section, we will first derive control laws that, in the absence of constraints,
will lead to overall system stability. Having derived these control laws, we will then show that there exist $\lambda$-
contractive sets $G_{1}$ and $G_{2}$, as described in \cite{lambda_contractive}, such that once $x_{1}^{aug}$ and $\tilde{x}$ enter
$G_{1},G_{2}$, they remain there (and in fact are driven further into the sets at the
next instant).

We consider two options for inner loop terminal control design, namely:
\begin{itemize}
  \item \emph{Exact reference model matching} - We design the inner loop control law such that $\tilde{x}(k+1)=A_{f}\tilde{x}(k)$, which guarantees that $v(k)$ asymptotically tracks $v_{des}^{f}(k)$.
  \item \emph{Approximate reference model matching} - The inner loop control law is designed such that the closed inner loop
  is stable and a small gain condition is satisfied.
\end{itemize}

\subsection{Terminal Control Law Design with Exact Reference Model Matching}
\begin{figure}[thpb]
\centering
\includegraphics[width=3.5in]{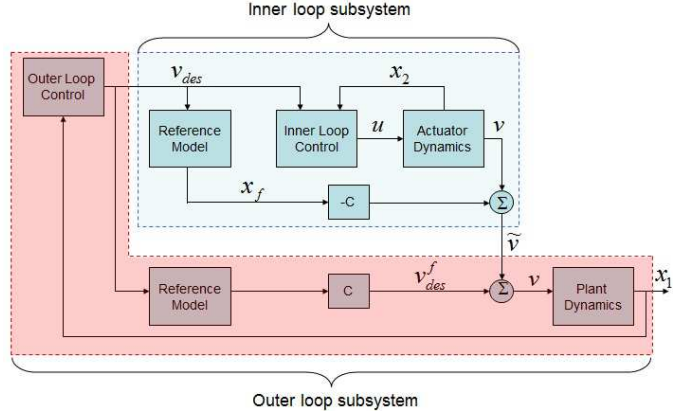}
\caption{Block diagram of the hierarchical control strategy under terminal control laws, rearranged for
analysis purposes.} \label{hierarchical_rm_analysis}
\end{figure}

In the case of exact reference model matching, in order to derive a model-matching controller, we assume that
the reference model is cast in a specific form that is compatible with the actuator dynamics; specifically, we assume that:
\begin{itemize}
  \item \emph{Assumption 6.} $A_{f}$ in (\ref{ref_model}) is written in the same block CCF as $A_{2}$ (as described in \cite{miso_ccf}).
\end{itemize}
\begin{itemize}
  \item \emph{Assumption 7.} Taking $R_{2}$ and $R_{f}$ as the set of rows of $B_{2}$ and $B_{f}$, respectively, that contain nonzero entries (which also correspond to the full rows of $A_{2}$ and $A_{f}$), we assume that $R_{f} \subset R_{2}$, i.e., each nonzero row of $B_{f}$ is also a nonzero row of $B_{2}$.
\end{itemize}
\noindent These assumptions, in conjunction with Assumptions 1-5, place restrictions on $A_{f}$ and $B_{f}$ that ensure that a stabilizing, reference model-matching
inner loop control law can be designed. In particular,it is
possible to design outer and inner loop terminal control laws with desirable properties, according to the following
proposition:
\begin{prop} \label{terminal_control} (Terminal control laws for exact matching): Given that Assumptions 1-7 hold, there exist control laws $v_{des}(k) = -K_{1}x_{1}^{aug}(k)$ and $u(k) = K_{21}v_{des}(k) - K_{22}x_{2}(k)$
which, when substituted into (\ref{ol_aug_dynamics}) and (\ref{il_error_dynamics}), yield:
\begin{eqnarray}\label{full_cl_dynamics}
\nonumber x_{1}^{aug}(k+1) &=& (A_{1}^{aug}-B_{1}^{aug}K_{1})x_{1}^{aug}(k) \\
&& + B_{1}^{aug}\tilde{v}(k), \\
\nonumber \tilde{x}(k+1) &=& A_{f}\tilde{x}(k),
\end{eqnarray}
\noindent where $\|\bar{\lambda}_{i}(A_{1}^{aug}-B_{1}^{aug}K_{1}))\| < 1, \forall i$,
and render the origin of the overall system, $x_{1}^{aug}=0,\tilde{x}=0$, asymptotically stable. \end{prop}
\begin{pf} Since the pair $(A_{1},B_{1})$ is controllable and the reference model does not share zeros with unstable poles of $A_{1}$, it follows that the pair $(A_{1}^{aug},B_{1}^{aug})$ is
stabilizable.  Thus, $K_{1}$ can be designed to ensure that $\|\bar{\lambda}_{i}(A_{1}^{aug}-B_{1}^{aug}K_{1}))\| < 1, \forall i$.

To show the second part of the proposition, recall that the inner loop dynamics are
expressed in (\ref{il_error_dynamics}) by:
\begin{equation}\label{inner_dynamics_again}
\tilde{x}(k+1) = A_{2}\tilde{x}(k)+(A_{2}-A_{f})x_{f}(k)+B_{2}u(k)-B_{f}v_{des}(k).
\end{equation}
It follows from the block CCF of $A_{2}$, $B_{2}$, in conjunction with Assumptions 6 and 7 (which impose a suitable block CCF structure on $A_{f}$ and $B_{f}$), that we can choose $K_{21}$ and $K_{22}$ to satisfy:
\begin{eqnarray}
B_{2}K_{21} &=& B_{f}, \\
\nonumber B_{2}K_{22} &=& A_{2} - A_{f},
\end{eqnarray}
\noindent which, when substituted into the inner loop dynamics, yields:
\begin{equation}
\tilde{x}(k+1) = A_{f} \tilde{x}(k).
\end{equation}
\noindent To see this, let $R_{2}$ be the indices corresponding to the nonzero rows of $B_{2}$,
and let $R_{f}$ be the indices corresponding to the nonzero rows of $B_{f}$. Because Assumption 7 requires that that $R_{f} \subseteq R_{2}$, it possible to achieve $B_{2}K_{21} = B_{f}$. Furthermore, let
$i$ represent any zero row of $B_{2}$ (and also $B_{f}$). It follows from block CCF (imposed by Assumption 6) and Assumption 7 that $A_{2ij} = A_{fij} \forall j$, making
it possible to achieve $B_{2}K_{22} = A_{2} - A_{f}$.

Because $\|\bar{\lambda}_{i}(A_{1}^{aug}-B_{1}^{aug}K_{1}))\| < 1, \forall i$ and
$\|\bar{\lambda}_{i}(A_{f}))\| < 1, \forall i$, it follows that both the closed inner
and outer loops (\ref{full_cl_dynamics}) are input-to-state stable (ISS) under the aforementioned control laws.  For small gain analysis, it is convenient to recast the system block diagram of Fig. \ref{hierarchical_basic} in the nonminimal representation of Fig. \ref{hierarchical_rm_analysis}, where offsetting copies of the reference model are embedded in both the inner and outer loops. Closed outer loop stability
guarantees a finite $l_{2}$ gain, $\gamma_{1}$, from $\tilde{v}$ to $v_{des}$, and exact
inner loop reference model matching guarantees an $l_{2}$ gain of $\gamma_{2}=0$, from
$v_{des}$ to $\tilde{v}$.  Therefore, the small gain condition, $\gamma_{1}\gamma_{2}<1$, is
satisfied, and in conjunction with outer and inner loop ISS, this proves asymptotic
stability of $x_{1}^{aug}=0,\tilde{x}=0$.
\begin{flushright} $\Box$
\end{flushright}
\end{pf}
\subsection{Terminal Control Law Design with Inexact Reference Model Matching and Small Gain Condition}
For many systems, such as non-minimum phase
systems and high-order, high relative degree systems, exact reference model matching is unrealistic.  Exact matching
is not essential, however, as shown in the following Proposition:
\begin{prop} \label{terminal_control_inexact} (Terminal control laws for inexact matching): Given that Assumptions 1-5 hold, there exist control laws
$v_{des}(k) = -K_{1}x_{1}^{aug}(k)$ and $u(k) = K_{21}v_{des}(k) - K_{22}x_{2}(k)$
which, when substituted into (\ref{ol_aug_dynamics}) and (\ref{il_error_dynamics}), yield:
\begin{eqnarray}\label{full_inexact_cl_dynamics}
\nonumber x_{1}^{aug}(k+1) &=& (A_{1}^{aug}-B_{1}^{aug}K_{1})x_{1}^{aug}(k) + B_{1}^{aug}\tilde{v}(k), \\
\tilde{x}(k+1) &=& (A_{2}-B_{2}K_{22})\tilde{x}(k) \\
\nonumber && + (A_{2}-A_{f}-B_{2}K_{22})x_{f}(k) \\
\nonumber && + (B_{2}K_{21}-B_{f})v_{des}(k),
\end{eqnarray}
\noindent where $\|\bar{\lambda}_{i}(A_{1}^{aug}-B_{1}^{aug}K_{1}))\| < 1, \forall i$ and
$\|\bar{\lambda}_{i}(A_{2}-B_{2}K_{22}))\| < 1, \forall i$.  Furthermore, suppose that
$K_{1}$, $K_{21}$, and $K_{22}$ are designed so that the $l_{2}$ gains from $\tilde{v}$ to
$v_{des}$ and $v_{des}$ to $\tilde{v}$, denoted $\gamma_{1}$ and $\gamma_{2}$, respectively,
satisfy the small gain condition, $\gamma_{1}\gamma_{2} < 1$.  Then the origin of the
overall system, $x_{1}^{aug}=0,\tilde{x}=0$, is asymptotically stable. \end{prop}
\begin{pf} Since the pair $(A_{1},B_{1})$ is controllable and the reference model does not share zeros with unstable poles of $A_{1}$, it follows that the pair $(A_{1}^{aug},B_{1}^{aug})$ is
stabilizable.  Thus, $K_{1}$ can be designed to ensure that $\|\bar{\lambda}_{i}(A_{1}^{aug}-B_{1}^{aug}K_{1}))\| < 0, \forall i$. From Assumption 2, the pair ($A_{2},B_{2}$) is controllable, and therefore there
exists $K_{2}$ such that $\|\bar{\lambda}_{i}(A_{2}-B_{2}K_{2}))\| < 0, \forall i$. Thus, both the inner and outer closed-loop
dynamics of (\ref{full_inexact_cl_dynamics}), are input-to-state stable (ISS). By the hypotheses of Proposition \ref{terminal_control_inexact}, the small gain condition, i.e., $\gamma_{1}\gamma_{2}<1$, is satisfied. Together with ISS, this
proves asymptotic stability of $x_{1}^{aug}=0,\tilde{x}=0$.
\begin{flushright} $\Box$
\end{flushright}
\end{pf}
Since $\|\bar{\lambda}_{i}(A_{1}^{aug}-B_{1}^{aug}K_{1}))\| < 1, \forall i$ and $\|\bar{\lambda}_{i}(A_{2}-B_{2}K_{22}))\| < 1, \forall i$ it follows (\cite{hierarchical1}) that there exist quadratic
Lyapunov functions, $V(x_{1}^{aug}) = x_{1}^{aug,T}Qx_{1}^{aug}$ and $V(\tilde{x}) = \tilde{x}^{T}P\tilde{x}$, where $Q$ and $P$ are
positive definite symmetric matrices, such that when $v_{des}(k)=-K_{1}x_{1}^{aug}(k)$ and $u(k) = K_{21}v_{des}(k) - K_{22}x_{2}(k)$:
\begin{eqnarray}\label{lyap_lambda_inexact}
V_{1}(x_{1}^{aug}(k+1))-V_{1}(x_{1}^{aug}(k)) &<& -\alpha_{1} V_{1}(x_{1}^{aug}(k)) \\
\nonumber && + \bar{\gamma}_{1} \|\tilde{v}(k)\|^{2},
\end{eqnarray}
\begin{eqnarray}\label{lyap_alpha_inexact}
V_{2}(\tilde{x}(k+1))-V_{2}(\tilde{x}(k)) &<& -\alpha_{2}V_{2}(\tilde{x}(k)) \\
\nonumber && + \bar{\gamma}_{21} \|v_{des}(k)\|^{2} + \bar{\gamma}_{22} \|x_{f}(k)\|^{2}.
\end{eqnarray}
\noindent for some $\alpha_{1} > 0$, $\alpha_{2} > 0$, $\bar{\gamma}_{1} > 0$, $\bar{\gamma}_{21} > 0$, and $\bar{\gamma}_{22} \geq 0$ ($\bar{\gamma}_{22} = 0$ under exact reference model matching). This fact will
be important in demonstrating the existence and construction of $\lambda$-contractive terminal
constraint sets.

\subsection{Design of Terminal Constraint Sets Under Exact Reference Model Matching}
Now, we show that $\lambda$-contractive sets, $G_{1}$ and $G_{2}$, conforming to the definition in \cite{lambda_contractive},
exist for the outer and inner loops. To guarantee that such sets exist, we make the following
trivial assumption regarding the feasible control input set, $U$:
\begin{itemize}
  \item \emph{Assumption 8.} $u=0$ lies in the interior of $U$.
\end{itemize}
\noindent We now demonstrate the existence of $\lambda$-contractive sets through
the following proposition:
\begin{prop} \label{lambda_contractive} (Existence of $\lambda$-contractive sets): Under Assumptions 1-8, there
exist sets $G_{1} \subset \mathbb{R}^{n_{1}+n_{2}}$ and $G_{2} \subset \mathbb{R}^{n_{2}}$, along
with scalars $\lambda_{1}: 0 \leq \lambda_{1} < 1$ and $\lambda_{2}: 0 \leq \lambda_{2} < 1$ such that if:
\begin{eqnarray}
\nonumber x_{1}^{aug}(k) &\in& G_{1}, \\
\tilde{x}(k) &\in& G_{2}, \\
\nonumber v_{des}(k) &=& -K_{1}x_{1}(k), \\
\nonumber u(k) &=& K_{21}v_{des}(k) - K_{22}x_{2}(k),
\end{eqnarray}
then:
\begin{eqnarray}
\nonumber u(k) &\in& U, \\
x_{1}^{aug}(k+1) &\in& \lambda_{1}G_{1}, \\
\nonumber \tilde{x}(k+1) &\in& \lambda_{2}G_{2}.
\end{eqnarray} \end{prop}
\begin{pf} To construct $G_{1}$, take:
\begin{equation}\label{G_1_def}
G_{1} \triangleq \{x_{1}^{aug}: V_{1}(x_{1}^{aug}) < V_{1}^{*}\},
\end{equation}
\noindent where $V_{1}^{*} > 0$. It follows from the continuity of $V_{1}(x_{1}^{aug})$ that there exists some
$\lambda_{1}: 0 \leq \lambda_{1} < 1$, $\lambda_{1}^{*}: 0 \leq \lambda_{1}^{*} < 1$, $\epsilon_{1} > 0$ such that:
\begin{eqnarray}\label{epsilon}
\lambda_{1}^{*} &>& 1-\alpha_{1}+\epsilon_{1}, \\
\nonumber V_{1}(x_{1}^{aug}(k+1))&\leq& \lambda_{1}^{*}V_{1}^{*} \Rightarrow x_{1}(k+1) \in \lambda_{1}G_{1}.
\end{eqnarray}
It follows from (\ref{lyap_lambda_inexact}), (\ref{G_1_def}), and (\ref{epsilon}) that if:
\begin{equation}\label{v_tilde_limits}
\|\tilde{v}(k)\|^{2} \leq \frac{\epsilon_{1} V_{1}^{*}}{\bar{\gamma}_{1}},
\end{equation}
and $x_{1}^{aug}(k) \in G_{1}$, then:
\begin{eqnarray}\label{last_G1}
V_{1}(x_{1}^{aug}(k+1)) &\leq& \lambda_{1}^{*} V_{1}^{*}, \\
\nonumber x_{1}^{aug}(k+1) &\in& \lambda_{1}G_{1}.
\end{eqnarray}
\noindent To see this, note that (\ref{lyap_lambda_inexact}) can be rearranged as:
\begin{eqnarray}
V_{1}(x_{1}^{aug}(k+1)) &<& (1-\alpha_{1}+\epsilon_{1}) V_{1}^{*} \\
\nonumber && + (1-\alpha_{1})(V_{1}(x_{1}^{aug}(k))-V_{1}^{*}),
\end{eqnarray}
\noindent and noting that $(1-\alpha_{1})(V_{1}(x_{1}^{aug}(k))-V_{1}^{*}) \leq 0$ when $x_{1}^{aug} \in G_{1}$,
it follows that:
\begin{equation}\label{last_G1}
V_{1}(x_{1}^{aug}(k+1)) < (1-\alpha_{1}+\epsilon_{1}) V_{1}^{*}.
\end{equation}

To construct $G_{2}$, take:
\begin{equation}
G_{2} \triangleq \{\tilde{x}: V_{2}(\tilde{x})\leq V_{2}^{*}\},
\end{equation}
\noindent where $V_{2}^{*} > 0$. It follows from (\ref{lyap_alpha_inexact}) and the
continuity of $V_{2}(\tilde{x})$ that if $\tilde{x}(k) \in G_{2}$ and
$u(k)=u_{t}(k)$, then:
\begin{eqnarray}
V_{2}(\tilde{x}(k+1)) &<& (1-\alpha_{2})V_{2}^{*}, \\
\nonumber \tilde{x}(k+1) &\in& \lambda_{2} G_{2},
\end{eqnarray}
\noindent for some $\lambda_{2}: 0 \leq \lambda_{2} < 1$.

It remains to select $V_{1}^{*}$ and $V_{2}^{*}$ such that
$u_{t}(k) \in U, \forall x_{1}^{aug}(k) \in G_{1}, \tilde{x}(k) \in G_{2}$, and
$\|\tilde{v}(k)\|^{2}$ satisfies (\ref{v_tilde_limits}) whenever $\tilde{x}(k) \in G_{2}$.

To ensure that $u_{t}(k) \in U$, note that the inner loop terminal control law (\ref{inner_terminal})
can be written as:
\begin{equation}\label{inner_control_rewritten}
u_{t}(k) = (-K_{21}K_{1} + K_{22})x_{f}(k) - K_{22}\tilde{x}(k),
\end{equation}
\noindent and that
\begin{equation}\label{inner_control_norms}
\|u_{t}(k)\| \leq (\|K_{21}K_{1}\|+\|K_{22}\|)\|x_{1}^{aug}(k)\| + \|K_{22}\|\|\tilde{x}(k)\|.
\end{equation}
It follows from (\ref{inner_control_norms}) and Assumption 8 that one can choose $x_{1}^{max} > 0$ and $\tilde{x}^{max} > 0$ such that whenever
$\|x_{1}^{aug}(k)\| \leq x_{1}^{max}$ and $\|\tilde{x}(k)\| \leq \tilde{x}^{max}$, $u_{t}(k) \in U$.
From the quadratic structure of $V_{1}(x_{1}^{aug})$, it follows that
$\|x_{1}^{aug}(k)\| \leq x_{1}^{max}$ whenever
\begin{equation}
V_{1}(x_{1}^{aug}(k)) \leq \bar{\lambda}_{min}(Q)(x_{1}^{max})^{2},
\end{equation}
\noindent where $\bar{\lambda}_{min}(Q)$ is the smallest eigenvalue of $Q$. Therefore, taking
\begin{equation}\label{V_1_star_constraint}
V_{1}^{*} \leq \bar{\lambda}_{min}(Q)(x_{1}^{max})^{2}
\end{equation}
\noindent guarantees that $\|x_{1}^{aug}\| \leq x_{1}^{max}$ whenever $x_{1}^{aug} \in G_{1}$. Similarly,
from the quadratic structure of $V_{2}(\tilde{x})$, it follows that
$\|\tilde{x}(k)\| \leq \tilde{x}^{max}$ whenever
\begin{equation}
V_{2}(\tilde{x}(k)) \leq \bar{\lambda}_{min}(P)(\tilde{x}^{max})^{2},
\end{equation}
\noindent where $\bar{\lambda}_{min}(P)$ is the smallest eigenvalue of $P$. Therefore, taking
\begin{equation}\label{V_2_star_constraint1}
V_{2}^{*} \leq \bar{\lambda}_{min}(P)(\tilde{x}^{max})^{2}
\end{equation}
\noindent guarantees that $\|\tilde{x}\| \leq \tilde{x}^{max}$ whenever $\tilde{x} \in G_{2}$.
(\ref{V_1_star_constraint}) and (\ref{V_2_star_constraint1}) together guarantee that $u_{t}(k) \in U$
whenever $x_{1}^{aug}(k) \in G_{1}, \tilde{x}(k) \in G_{2}$.

Finally, $V_{2}^{*}$ needs to be selected so that $\|\tilde{v}(k)\|^{2}$ satisfies (\ref{v_tilde_limits}) whenever $\tilde{x}(k) \in G_{2}$. Manipulation of (\ref{v_tilde_limits}) shows that this is the case when
\begin{equation}\label{first_v_tilde_eq}
\|\tilde{x}(k)\|^{2} \leq \frac{\epsilon_{1} V_{1}^{*}}{\|C\|^{2}\bar{\gamma}_{1}},
\end{equation}
\noindent and it follows from the quadratic structure of $V_{2}(\tilde{x})$ that (\ref{first_v_tilde_eq}) is
satisfied whenever
\begin{equation}
V_{2}(\tilde{x}(k)) \leq \frac{\epsilon_{1} V_{1}^{*}\bar{\lambda}_{min}(P)}{\|C\|^{2}\bar{\gamma}_{1}}.
\end{equation}
\noindent Substituting (\ref{V_1_star_constraint}) for $V_{1}^{*}$, it follows that
by taking
\begin{equation}
V_{2}^{*} \leq \frac{\epsilon_{1} (x_{1}^{max})^{2}\bar{\lambda}_{min}(Q)\bar{\lambda}_{min}(P)}{\bar{\gamma}_{1} \|C\|^{2}},
\end{equation}
\noindent one guarantees that (\ref{v_tilde_limits}) is satisfied whenever $\tilde{x}(k) \in G_{2}$. In
order to simultaneously ensure that $u_{t}(k) \in U$, we take:
\begin{equation}\label{V_2_star_constraint}
V_{2}^{*} = \min \{(\tilde{x}^{max})^{2}\bar{\lambda}_{min}(P),\frac{\epsilon_{1} (x_{1}^{max})^{2}\bar{\lambda}_{min}(Q)\bar{\lambda}_{min}(P)}{\bar{\gamma}_{1} \|C\|^{2}}\}.
\end{equation}
\begin{flushright} $\Box$
\end{flushright}
\end{pf}
The proof of Proposition \ref{lambda_contractive} is constructive in the sense that it provides the method by which one can construct $G_{1}$ and $G_{2}$, and determine suitable values for $\lambda_{1}$ and $\lambda_{2}$, respectively.

\subsection{Design of Terminal Constraint Sets Under Inexact Reference Model Matching}
It is also possible to derive constraint sets $G_{1}$ and $G_{2}$ under inexact, but sufficiently accurate reference model matching. The existence of $\lambda$-contractive constraint sets and the conditions under which
they are guaranteed to exist are given in the following proposition:
\begin{prop} \label{lambda_contractive_inexact} (Existence of $\lambda$-contractive sets with inexact reference
model matching): Suppose that Assumptions 1-5 and Assumption 8 hold. Furthermore, suppose that $V_{1}(x_{1}^{aug})$ and $V_{2}(\tilde{x})$, along with scalars $\alpha_{1}$, $\alpha_{2}$, $\bar{\gamma}_{1}$, $\bar{\gamma}_{21}$, and $\bar{\gamma}_{22}$
from (\ref{lyap_lambda_inexact}) and (\ref{lyap_alpha_inexact}), along with matrices $C$ and $K_{1}$, satisfy the following inequality:
\begin{equation}\label{small_gain_constraint_inequality}
\alpha_{1}\alpha_{2}\bar{\lambda}_{min}(P)\bar{\lambda}_{min}(Q) \geq \bar{\gamma}_{1}\|C\|^{2}(\bar{\gamma}_{21}\|K_{1}\|^{2}+\bar{\gamma}_{22}),
\end{equation}
\noindent where $\bar{\lambda}_{min}(P)$ and $\bar{\lambda}_{min}(Q)$ are the minimum
eigenvalues of $P$ and $Q$, respectively. Then there exist sets $G_{1} \subset \mathbb{R}^{n_{1}+n_{2}}$ and $G_{2} \subset \mathbb{R}^{n_{2}}$, along
with scalars $\lambda_{1}: 0 \leq \lambda_{1} < 1$ and $\lambda_{2}: 0 \leq \lambda_{2} < 1$ such that if:
\begin{eqnarray}
\nonumber x_{1}^{aug}(k) &\in& G_{1}, \\
\tilde{x}(k) &\in& G_{2}, \\
\nonumber v_{des}(k) &=& -K_{1}x_{1}(k), \\
\nonumber u(k) &=& K_{21}v_{des}(k) - K_{22}x_{2}(k),
\end{eqnarray}
then:
\begin{eqnarray}
\nonumber u(k) &\in& U, \\
x_{1}^{aug}(k+1) &\in& \lambda_{1}G_{1}, \\
\nonumber \tilde{x}(k+1) &\in& \lambda_{2}G_{2}.
\end{eqnarray} \end{prop}
\begin{pf} The construction of $G_{1}$ is done identically
to Proposition \ref{lambda_contractive}, taking:
\begin{equation}\label{G_1_def_inexact}
G_{1} \triangleq \{x_{1}^{aug}: V_{1}(x_{1}^{aug}) < V_{1}^{*}\},
\end{equation}
\noindent where $V_{1}^{*} > 0$. Equations (\ref{epsilon})-(\ref{last_G1}) remain unchanged and follow the
same derivation as in Proposition \ref{lambda_contractive}.

For the construction of $G_{2}$, we take:
\begin{equation}\label{G_2_def_inexact}
G_{2} \triangleq \{\tilde{x}: V_{2}(\tilde{x})\leq V_{2}^{*}\},
\end{equation}
\noindent where $V_{2}^{*} > 0$. It follows from the continuity of $V_{2}(\tilde{x})$ that there exists some
$\lambda_{2}: 0 \leq \lambda_{2} < 1$, $\lambda_{2}^{*}: 0 \leq \lambda_{2}^{*} < 1$, $\epsilon_{2} > 0$ such that:
\begin{eqnarray}\label{epsilon_inexact}
\lambda_{2}^{*} &>& 1-\alpha_{2}+\epsilon_{2}, \\
\nonumber V_{2}(\tilde{x}(k+1))&\leq& \lambda_{2}^{*}V_{2}^{*} \Rightarrow \tilde{x}(k+1) \in \lambda_{2}G_{2}.
\end{eqnarray}
It follows from (\ref{lyap_alpha_inexact}), (\ref{G_2_def_inexact}), and (\ref{epsilon_inexact}) that if:
\begin{equation}\label{v_des_limits}
\bar{\gamma}_{21}\|v_{des}(k)\|^{2}+\bar{\gamma}_{22}\|x_{f}(k)\|^{2} \leq \epsilon_{2}V_{2}^{*},
\end{equation}
and $\tilde{x}(k) \in G_{2}$, then:
\begin{eqnarray}
V_{2}(\tilde{x}(k+1)) &\leq& \lambda_{2}^{*} V_{2}^{*}, \\
\nonumber \tilde{x}(k+1) &\in& \lambda_{2}G_{2}.
\end{eqnarray}

It remains to select $V_{1}^{*}$ and $V_{2}^{*}$ such that
$u_{t}(k) \in U, \forall x_{1}^{aug}(k) \in G_{1}, \tilde{x}(k) \in G_{2}$, and
$\|\tilde{v}(k)\|^{2}$ satisfies (\ref{v_tilde_limits}) whenever $\tilde{x}(k) \in G_{2}$. This
derivation is exactly the same here as in Proposition \ref{lambda_contractive}, and
(\ref{inner_control_rewritten})-(\ref{V_2_star_constraint1}) all hold.

Finally, $V_{1}^{*}$ needs to be selected so that $\|v_{des}(k)\|^{2}$ and $\|x_{f}(k)\|^{2}$ satisfy (\ref{v_des_limits}) whenever
$x_{1}^{aug} \in G_{1}$,and $V_{2}^{*}$ needs to be selected so that $\|\tilde{v}(k)\|^{2}$ satisfies (\ref{v_tilde_limits}) whenever $\tilde{x}(k) \in G_{2}$. For $V_{2}^{*}$, the derivation is the same as in Proposition \ref{lambda_contractive} and the
requirement is given by:
\begin{equation}\label{V2_inexact}
V_{2}^{*} \leq \frac{\epsilon_{1} V_{1}^{*}\bar{\lambda}_{min}(P)}{\|C\|^{2}\bar{\gamma}_{1}}.
\end{equation}
\noindent For $V_{1}^{*}$, we begin by noting that if:
\begin{equation}\label{V1_inexact}
V_{1}^{*} \leq \frac{\epsilon_{2} V_{2}^{*}\bar{\lambda}_{min}(Q)}{\|K_{1}\|^{2}\bar{\gamma}_{21}+\bar{\gamma}_{22}},
\end{equation}
\noindent then (\ref{v_des_limits}) is satisfied. To see this, note first that whenever $x_{1}^{aug}(k) \in G_{1}$, it follows from the quadratic form of $V_{1}(x_{1}^{aug})$ that:
\begin{equation}
\bar{\lambda}_{min}(Q)\|x_{1}^{aug}(k)\|^{2} \leq V_{1}^{*},
\end{equation}
\noindent from which it follows from substitution into (\ref{V1_inexact}) that:
\begin{equation}
\|x_{1}^{aug}(k)\|^{2} \leq \frac{\epsilon_{2}V_{2}^{*}}{\|K_{1}\|^{2}\bar{\gamma}_{21}+\bar{\gamma}_{22}}.
\end{equation}
\noindent Noting that $\bar{\gamma}_{21}\|v_{des}(k)\|^{2}+\bar{\gamma}_{22}\|x_{f}(k)\|^{2} \leq \|x_{1}^{aug}(k)\|^{2}(\|K_{1}\|^{2}\bar{\gamma}_{21}+\bar{\gamma}_{22})$, we can see immediately that (\ref{v_des_limits}) is satisfied.

Combining (\ref{V2_inexact}) and (\ref{V1_inexact}) with the requirements of
(\ref{V_1_star_constraint}) and (\ref{V_2_star_constraint}) gives the following two nonlinear equations that must be solved
for $V_{1}^{*}$ and $V_{2}^{*}$:
\begin{equation}\label{V_1_star_constraint_inexact}
V_{1}^{*} = \min \{(x_{1}^{max})^{2}\bar{\lambda}_{min}(Q),\frac{\epsilon_{2} V_{2}^{*}\bar{\lambda}_{min}(Q)}{\|K_{1}\|^{2}\bar{\gamma}_{21}+\bar{\gamma}_{22}}\}.
\end{equation}
\begin{equation}\label{V_2_star_constraint_inexact}
V_{2}^{*} = \min \{(\tilde{x}^{max})^{2}\bar{\lambda}_{min}(P),\frac{\epsilon_{1} V_{1}^{*}\bar{\lambda}_{min}(P)}{\|C\|^{2}\bar{\gamma}_{1}}\}.
\end{equation}
\noindent (\ref{V_1_star_constraint_inexact}) and (\ref{V_2_star_constraint_inexact}) will only admit
a solution if:
\begin{equation}\label{inexact_constraint_epsilon}
\frac{\epsilon_{1} \bar{\lambda}_{min}(P)}{\|C\|^{2}\bar{\gamma}_{1}} \geq \frac{ \|K_{1}\|^{2}\bar{\gamma}_{21}+\bar{\gamma}_{22}}{\epsilon_{2}\bar{\lambda}_{min}(Q)}.
\end{equation}
\noindent Noting that the only requirements on $\epsilon_{1}$ and $\epsilon_{2}$ are that
$\epsilon_{1}<\alpha_{1}$ and $\epsilon_{2}<\alpha_{2}$, $\epsilon_{1}$ and $\epsilon_{2}$ in
(\ref{inexact_constraint_epsilon}) can be replaced with $\alpha_{1}$ and $\alpha_{2}$, and
(\ref{inexact_constraint_epsilon}) can be rearranged to yield the constraint:
\begin{equation}\label{small_gain_constraint_inequality_repeat}
\alpha_{1}\alpha_{2}\bar{\lambda}_{min}(P)\bar{\lambda}_{min}(Q) \geq \bar{\gamma}_{1}\|C\|^{2}(\bar{\gamma}_{21}\|K_{1}\|^{2}+\bar{\gamma}_{22}),
\end{equation}
\noindent completing the proof.
\begin{flushright} $\Box$
\end{flushright}
\end{pf}
\noindent The proof of Proposition \ref{lambda_contractive_inexact} follows similar arguments to that of
Proposition \ref{lambda_contractive}, with the exception that now $V_{1}^{*}$ and $V_{2}^{*}$, which define the boundaries of $G_{1}$ and $G_{2}$, must satisfy
two coupled equations, and a solution to these coupled equations only exists when (\ref{small_gain_constraint_inequality})
is satisfied. Qualitatively speaking, satisfaction of (\ref{small_gain_constraint_inequality}) depends on two factors:
\begin{enumerate}
  \item Free response speed of the outer and inner loop systems, indicated by $\alpha_{1}$ and $\alpha_{2}$;
  \item Level of coupling between the outer and inner loop systems, indicated by
      $\bar{\gamma}_{1}$, $\bar{\gamma}_{21}$, and $\bar{\gamma}_{22}$.
\end{enumerate}

\section{Deriving Rate-Like Constraints on Control Inputs and Desired Virtual Control Inputs}\label{rate_constraints}
The results of Section \ref{lambda_contractive} provide a means by which outer
and inner loop control laws can be designed to yield local stability of the
origin of the overall system, i.e., $x_{1}^{aug}=0$, $\tilde{x}=0$. The
MPC optimizations of (\ref{mpc_outer_opt})-(\ref{mpc_outer_cost}) and (\ref{mpc_inner_opt})-(\ref{mpc_inner_cost})
are employed in order to expand the region of attraction beyond the intersection of
$G_{1}$ and $G_{2}$. In order to guarantee convergence to $G_{1}$ and $G_{2}$, the MPC
optimizations must not only impose a terminal constraint but must also ensure that
optimized trajectories do not differ too much from one time step to the next in order to
ultimately guarantee persistent feasibility of the optimization. This assurance is accomplished
through the imposition of rate-like constraints presented in this section. These rate-like constraints,
$\delta_{v_{des}}^{max}$ and $\delta_{u}^{max}$, which limit the variation of $\mathbf{v_{des}}$ and
$\mathbf{u}$ trajectories from one time instant to the next.
%

We begin with the following proposition, which follows from examination of the time series representation of the $x_{1}^{aug}$ trajectory:
\begin{prop} \label{OL_robustness_to_v_tilde_variation} (Robustness of outer loop MPC to variation in $\mathbf{\tilde{v}}$): Suppose
that, given
\begin{equation}
\nonumber \mathbf{\tilde{v}}(k-1) = \left[\begin{array}{ccc} \mathbf{\tilde{v}}(k-1|k-1) & \ldots & \mathbf{\tilde{v}}(k+N-1|k-1)\end{array}\right],
\end{equation}
\noindent a trajectory
\begin{equation}
\nonumber \mathbf{v_{des}}(k) = \left[\begin{array}{ccc} \mathbf{v_{des}}(k|k) & \ldots & \mathbf{v_{des}}(k+N-1|k)\end{array}\right],
\end{equation}
\noindent is computed that yields $\mathbf{x_{1}^{aug}}(k+N|k) \in \lambda_{1}G_{1}$.  Then there exists $\epsilon_{\tilde{v}}^{max} > 0$ such that if $\|\mathbf{\tilde{v}}(k+i|k)-\mathbf{\tilde{v}}(k+i|k-1)\| \leq \epsilon_{\tilde{v}^{max}}, i=1 \ldots N-1$ and
$\mathbf{v_{des}}(k+i|k+1)= \mathbf{v_{des}}(k+i|k), i=1 \ldots N-1$, then $\mathbf{x_{1}^{aug}}(k+N|k+1) \in G_{1}$. \end{prop}
\begin{pf} At step $k$ the outer loop dynamics over the MPC horizon can be expressed as:
\begin{eqnarray}
\mathbf{x_{1}^{aug}}(k+i|k) &=& (A_{1}^{aug})^{i}x_{1}^{aug}(k) \\
\nonumber && +\sum_{j=0}^{i-1}(A_{1}^{aug})^{j}(B_{f}^{aug}\mathbf{v_{des}}(k+i-j-1|k)\\
\nonumber && +B_{1}^{aug}\mathbf{\tilde{v}}(k+i-j-1|k-1))
\end{eqnarray}
\noindent for $i=1 \ldots N$.  An analogous expression exists at step $k+1$.  When $\mathbf{v_{des}}(k+i|k)= \mathbf{v_{des}}(k+i|k+1), i=1 \ldots N-1$, the difference between the predicted
trajectories at steps $k$ and $k+1$ is then given by:
\begin{eqnarray}
\nonumber \mathbf{x_{1}^{aug}}(k+i|k+1)-\mathbf{x_{1}^{aug}}(k+i|k) &=& \sum_{j=0}^{i-1}(A_{1}^{aug})^{j}B_{1}^{aug}(\mathbf{\tilde{v}}(k+i-j-1|k) \\
\nonumber &&-\mathbf{\tilde{v}}(k+i-j-1|k-1)),
\end{eqnarray}
\noindent which, in the case that $\|\mathbf{\tilde{v}}(k+i|k)-\mathbf{\tilde{v}}(k+i|k-1)\| \leq \epsilon_{\tilde{v}^{max}}, i=1 \ldots N-1$, leads to the inequality:
\begin{equation}
\|\mathbf{x_{1}^{aug}}(k+i|k)-\mathbf{x_{1}^{aug}}(k+i|k+1)\| \leq \epsilon_{\tilde{v}}^{max}\sum_{j=0}^{i-1}\|(A_{1}^{aug})^{j}B_{1}^{aug}\|.
\end{equation}
\noindent Since $\lambda_{1} < 1$, it follows that there exists $\Delta_{1} > 0$
such that if $\mathbf{x_{1}^{aug}}(k+N|k) \in \lambda_{1}G_{1}$ and $\mathbf{x_{1}^{aug}}(k+N|k)-\mathbf{x_{1}^{aug}}(k+N|k+1) \leq \Delta_{1}$, then
$\mathbf{x_{1}^{aug}}(k+N|k+1) \in G_{1}$.  Thus, by taking:
\begin{equation}
\epsilon_{\tilde{v}}^{max} \leq \frac{\Delta_{1}}{\sum_{j=0}^{N-1}\|(A_{1}^{aug})^{j}B_{1}^{aug}\|}
\end{equation}
\noindent we guarantee that $\mathbf{x_{1}^{aug}}(k+N|k) \in G_{1}$.
\begin{flushright} $\Box$
\end{flushright}
\end{pf}
\noindent The proof relies on a time series representation of the $x_{1}^{aug}$ trajectory, which
demonstrates that the step-to-step variation in $\mathbf{x_{1}^{aug}}$ can be upper bounded by
restricting the variation in $\mathbf{\tilde{v}}$.

We arrive at a very similar conclusion regarding the robustness of the inner loop MPC to variation in $\mathbf{x_{f}}$:
\begin{prop} \label{IL_robustness_to_x_f_variation} (Robustness of inner loop MPC to variation in $\mathbf{x_{f}}$): Suppose that, given
\begin{equation}
\nonumber \mathbf{x_{f}}(k) = \left[\begin{array}{ccc} \mathbf{x_{f}}(k|k) & \ldots & \mathbf{x_{f}}(k+N|k)\end{array}\right],
\end{equation}
\noindent a trajectory
\begin{equation}
\nonumber \mathbf{u}(k) = \left[\begin{array}{ccc} \mathbf{u}(k|k) & \ldots & \mathbf{u}(k+N-1|k)\end{array}\right],
\end{equation}
\noindent is computed that yields $\mathbf{\tilde{x}}(k+N|k) \in \lambda_{2}G_{2}$.  Then there exists $\epsilon_{x_{f}}^{max} > 0$ such that if $\|\mathbf{x_{f}}(k+N|k+1)-\mathbf{x_{f}}(k+N|k)\| \leq \epsilon_{x_{f}}^{max}$ and $\mathbf{u}(k+i|k+1)= \mathbf{u}(k+i|k), i=1 \ldots N-1$, then $\mathbf{\tilde{x}}(k+N|k+1) \in G_{2}$. \end{prop}
\begin{pf} Taking $\mathbf{u}(k+i|k+1)=\mathbf{u}(k+i|k), i=1 \ldots N-1$ yields
$\mathbf{x_{2}}(k+i|k+1)=\mathbf{x_{2}}(k+i|k), i=1 \ldots N$. Thus,
\begin{equation} \label{x2_traj}
\mathbf{\tilde{x}}(k+i|k+1)-\mathbf{\tilde{x}}(k+i|k) = \mathbf{x_{f}}(k+i|k+1)-\mathbf{x_{f}}(k+i|k), i=1 \ldots N.
\end{equation}
\noindent Since $\lambda_{2}<1$, there exists $\Delta_{2}>0$ such that if
$\mathbf{\tilde{x}}(k+N|k) \in \lambda_{2}G_{2}$ and $\|\mathbf{\tilde{x}}(k+N|k+1)-\mathbf{\tilde{x}}(k+N|k)\| \leq \Delta_{2}$,
then
\begin{equation} \label{x_tilde_inthere}
\mathbf{\tilde{x}}(k+N|k+1) \in G_{2}.
\end{equation}
\noindent From (\ref{x2_traj}) and (\ref{x_tilde_inthere}) it follows that by
taking $\epsilon_{x_{f}}^{max} = \Delta_{2}$, we guarantee that $\mathbf{\tilde{x}}(k+N|k+1) \in G_{2}$.
\begin{flushright} $\Box$
\end{flushright}
\end{pf}
It is possible to convert the state constraints of Propositions \ref{OL_robustness_to_v_tilde_variation} and \ref{IL_robustness_to_x_f_variation}
to input constraints (on $\mathbf{v_{des}}$ and $\mathbf{u}$), which are easily enforced and will always result in a feasible
optimization problem (as opposed to state constraints, which are not in general guaranteed to result in a feasible constrained optimization).  These input constraints are given in the following propositions:
\begin{prop} \label{OL_robustness_to_u_variation} (Converting constraints on $\mathbf{\tilde{v}}$ to constraints
on $\mathbf{u}$): There exists $\delta_{u}^{max} > 0$ such that
if $\|\mathbf{u}(k+i|k)-\mathbf{u^{o}}(k+i|k-1)\|\leq \delta_{u}^{max}$, $i=0 \ldots N-2$, then $\|\mathbf{\tilde{v}}(k+i|k)-\mathbf{\tilde{v}}(k+i|k-1)\| \leq \epsilon_{\tilde{v}}^{max}$, $i=0 \ldots N-1$. \end{prop}
\begin{pf} For this proof, it is convenient to express the inner loop dynamics as:
\begin{equation}
\tilde{v}(k+1) = C(A_{2}x_{2}(k)+B_{2}u(k)-x_{f}(k+1)),
\end{equation}
\noindent from which it from a time series expansion that:
\begin{eqnarray}
\nonumber \mathbf{\tilde{v}}(k+i|k)-\mathbf{\tilde{v}}(k+i|k-1) &=& C\sum_{j=0}^{i-1}(A_{2}^{j}B_{2}(\mathbf{u}(k+i-j-1|k)\\
&& -\mathbf{u^{o}}(k+i-j-1|k-1))) \\
&& \nonumber -C(\mathbf{x_{f}}(k+i|k)-\mathbf{x_{f}}(k+i|k-1)),
\end{eqnarray}
\noindent and
\begin{equation}
\|\mathbf{\tilde{v}}(k+i|k)-\mathbf{\tilde{v}}(k+i|k-1)\| \leq \delta_{u}^{max} \|C\| \sum_{j=0}^{i-1}\|A_{2}^{j}B_{2}\| + \|C\| \epsilon_{x_{f}}^{max}.
\end{equation}
\noindent It follows that if we take:
\begin{equation}
\delta_{u}^{max} \leq \frac{\epsilon_{\tilde{v}^{max}}-\|C\| \epsilon_{x_{f}^{max}}}{\|C\| \sum_{j=0}^{N-1}\|A_{2}^{j}B_{2}\|},
\end{equation}
\noindent then we have $\|\mathbf{\tilde{v}}(k+i|k)-\mathbf{\tilde{v}}(k+i|k-1)\| \leq \epsilon_{\tilde{v}}^{max}, i = 0 \ldots N-1$.
\begin{flushright} $\Box$
\end{flushright}
\end{pf}
\noindent The proof uses the time series expression of the inner loop dynamics to demonstrate that one can restrict
the step-to-step variation in $\mathbf{u}$ and achieve the required bound on the step-to-step
variation in $\mathbf{\tilde{v}}$.

Constraints on $\mathbf{x_{f}}$ can similarly be converted to constraints on $\mathbf{v_{des}}$, as
presented in the following proposition:
\begin{prop} \label{IL_robustness_to_v_des_variation} (Converting constraints on $\mathbf{x_{f}}$ to constraints
on $\mathbf{v_{des}}$): There exists $\delta_{v_{des}}^{max} > 0$ such that
if $\|\mathbf{v_{des}}(k+i|k+1)-\mathbf{v_{des}}(k+i|k)\| \leq \delta_{v_{des}}^{max}$, $i=1 \ldots N-1$, then $\|\mathbf{x_{f}}(k+N|k+1)-\mathbf{x_{f}}(k+N|k)\| \leq \epsilon_{x_{f}}^{max}$. \end{prop} 
\begin{pf} Recall that
the reference model dynamics are given by:
\begin{equation}
x_{f}(k+1) = A_{f}x_{f}(k) + B_{f}v_{des}(k),
\end{equation}
\noindent from which it follows that:
\begin{eqnarray}
\nonumber \mathbf{x_{f}}(k+i|k+1)-\mathbf{x_{f}}(k+i|k) &=& \sum_{j=0}^{i-1}(A_{f}^{j}B_{f}(\mathbf{v_{des}}(k+i-j-1|k+1)\\
&& -\mathbf{v_{des}}(k+i-j-1|k))),
\end{eqnarray}
\noindent and
\begin{equation}
\|\mathbf{x_{f}}(k+i|k+1)-\mathbf{x_{f}}(k+i|k)\| \leq \delta_{v_{des}}^{max} \sum_{j=0}^{i-1}\|A_{f}^{j}B_{f}\|.
\end{equation}
\noindent If we take:
\begin{equation}
\delta_{v_{des}}^{max} \leq \frac{\epsilon_{x_{f}^{max}}}{\sum_{j=0}^{N-1}\|A_{f}^{j}B_{f}\|},
\end{equation}
\noindent then we have $\|\mathbf{x_{f}}(k+N|k+1)-\mathbf{x_{f}}(k+N|k)\| \leq \epsilon_{x_{f}}^{max}$.
\begin{flushright} $\Box$
\end{flushright}\end{pf}

\section{Persistent Feasibility, Convergence, and Stability}\label{feasibility_stability}
In this section, we show how the constraints derived in Sections 4 and 5
result in persistent feasibility of the MPC optimization problem and asymptotic
stability of the overall system, with a region of attraction that is identical to
the set of states for which the initial optimization problem is feasible.

\subsection{Persistent Feasibility}
Because the rate-like constraints cannot be applied
at step $k=0$ (since there is no step $k=-1$ against which to compare), we
make the following initial feasibility Assumption for step $k=0$:

\emph{Initial Feasibility Assumption}: There exists a set
$X \in \mathbb{R}^{n_{1}+2n_{2}}$, such that if $\left[\begin{array}{cc}x_{1}^{aug}(0)^{T} & \tilde{x}(0)^{T}\end{array}\right]^{T} \in X$, then
$\mathbf{v_{des}}(0)$ and $\mathbf{u}(0)$ can be chosen and are chosen such that
$|\mathbf{\tilde{v}}(i|0)-\mathbf{\tilde{v}}(i|-1)| \leq \epsilon_{\tilde{v}}, i=0 \ldots N-1$ and the MPC optimization problem is feasible.
%

Given this assumption, we now state the persistent feasibility result.
\begin{prop} \label{feasibility} (Persistent feasibility): Suppose that the initial conditions satisfy
$\left[\begin{array}{cc}x_{1}^{aug}(0)^{T} & \tilde{x}(0)^{T}\end{array}\right]^{T} \in X$.
Then both the outer and inner loop MPC optimizations are feasible at every step, $k \geq 0$. \end{prop}
\begin{pf} Feasibility at $k=0$ is guaranteed by the initial feasibility assumption.

\emph{Outer loop MPC feasibility for $k \geq 1$}: By inner loop constraint (\ref{mpc_inner_constraints}), combined
with Proposition \ref{OL_robustness_to_u_variation}, we guarantee that
$\|\mathbf{\tilde{v}}(k+i|k)-\mathbf{\tilde{v}}(k+i|k-1)\| \leq \epsilon_{\tilde{v}}^{max}$
for $i= 0 \ldots N-2$.  Thus, if we take $\mathbf{v_{des}}(k+i|k)=\mathbf{v_{des}}(k+i|k-1)$ for
$i=0 \ldots N-2$, then we achieve:
\begin{eqnarray}
\mathbf{x_{1}^{aug}}(k+N-1|k) & \in & G_{1}, \\
\nonumber \|\mathbf{v_{des}}(k+i|k)-\mathbf{v_{des}}(k+i|k-1)\| = 0 &\leq& \delta_{v_{des}}^{max} \beta^{k}, \\ \nonumber && i=0 \ldots N-2.
\end{eqnarray}
By construction of $G_{1}$ and $G_{2}$, taking $\mathbf{v_{des}^{o}}(k+N-1|k)=-K_{1}\mathbf{x_{1}^{aug}}(k+N-1|k)$
results in $\mathbf{x_{1}^{aug}}(k+N|k) \in \lambda G_{1}$.

\emph{Inner loop MPC feasibility for $k \geq 1$}: By outer loop constraint (\ref{mpc_outer_constraints}), combined
with Proposition \ref{IL_robustness_to_v_des_variation}, we guarantee that
$\|\mathbf{x_{f}}(k+i|k)-\mathbf{x_{f}}(k+i|k-1)\| \leq \epsilon_{x_{f}}^{max}$
for $i=0 \ldots N-2$.  Thus, if we take $\mathbf{u^{o}}(k+i|k)=\mathbf{u^{o}}(k+i|k-1)$ for
$i=0 \ldots N-2$, then we achieve:
\begin{eqnarray}
\nonumber \mathbf{\tilde{x}}(k+N-1|k) & \in & G_{2}, \\
\mathbf{u}(k+i|k) & \in & U, \\
\nonumber \|\mathbf{u}(k+i|k)-\mathbf{u^{o}}(k+i|k-1)\| = 0 &\leq& \delta_{u}^{max} \beta^{k}, \\
\nonumber && i=0 \ldots N-2.
\end{eqnarray}
Given that $\mathbf{x_{1}^{aug}}(k+N-1|k) \in G_{1}$, applying
$\mathbf{u}(k+N-1|k) = K_{21}\mathbf{v_{des}}(k+N-1|k) - K_{22}\mathbf{x_{f}}(k+N-1|k) - K_{22}\mathbf{\tilde{v}}(k+N-1|k)$ yields $\mathbf{\tilde{x}}(k+N|k) \in \lambda_{2}G_{2}$.
\begin{flushright} $\Box$
\end{flushright}
\end{pf}
\noindent The proof follows from the rate-like constraints imposed
on $\mathbf{v_{des}}(k)$ and $\mathbf{u}(k)$.  Specifically, if the variations in $\mathbf{v_{des}}$ and
$\mathbf{u}$ are sufficiently small from step $k$ to $k+1$, then the optimization
problem remains feasible at step $k+1$.

\subsection{Convergence}
Having shown that the optimization problems are persistently feasible, the next step is to
show that the control laws do in fact result in finite-time convergence to $G_{1}$ and $G_{2}$.  This is
given in the following proposition:
\begin{prop} \label{convergence} (Convergence to $G_{1}$, $G_{2}$): Suppose that the initial conditions satisfy
$\left[\begin{array}{cc}x_{1}^{aug}(0)^{T} & \tilde{x}(0)^{T}\end{array}\right] \in X$.
Then there exists a scalar integer $N^{*} > 0$ such that, after applying the MPC algorithm for $N^{*}$
steps, we have $x_{1}^{aug}(N^{*}) \in G_{1}$ and $\tilde{x}(N^{*}) \in G_{2}$. \end{prop}
\begin{pf} By the inner and outer loop rate-like constraints, we have:
\begin{eqnarray}
\|\mathbf{u^{o}}(k+i|k)-u(k+i)\| &\leq& i \delta_{u}^{max} \beta^{k}, \\
\nonumber \|\mathbf{v_{des}^{o}}(k+i|k)-v_{des}(k+i)\| &\leq& i \delta_{v_{des}}^{max} \beta^{k},
\end{eqnarray}
\noindent For the outer loop, it follows that:
\begin{eqnarray}
\nonumber \|\mathbf{x_{1}^{aug}}(k+N|k)-x_{1}^{aug}(k+N)\| &\leq& N(\sum_{j=0}^{i-1} \|(A_{1}^{aug})^{j}B_{1}^{aug}\| \epsilon_{\tilde{v}}^{max} \\
\nonumber && + \sum_{j=0}^{i-1} \|A_{1}^{j}B_{2}\| \epsilon_{x_{f}}^{max})\beta^{k},
\end{eqnarray}
\noindent which, after collecting constant terms into one lumped constant, $Q$, can be rewritten
compactly as:
\begin{equation}\label{ol_convergence_1}
\|\mathbf{x_{1}^{aug}}(k+N|k)-x_{1}^{aug}(k+N)\| \leq QN \beta^{k}.
\end{equation}
Because $\lambda_{1}G_{1} \in G_{1}$, there exists a positive scalar $\Delta x_{1}^{aug}$
such that for any two vectors $x_{1a}^{aug} \in \lambda_{1}G_{1}$ and $x_{1b}^{aug} \in G_{1}$,
$\|x_{1a}^{aug}-x_{1b}^{aug}\|<\Delta x_{1}^{aug}$. To guarantee that $\mathbf{x_{1}^{aug}}(k+N|k) \in \lambda_{1}G_{1} \Rightarrow x_{1}^{aug}(k+N) \in G_{1}$, it suffices to ensure that:
\begin{equation}\label{ol_convergence_2}
\|\mathbf{x_{1}^{aug}}(k+N|k)-x_{1}^{aug}(k+N)\| < \Delta x_{1}^{aug}.
\end{equation}

It follows through manipulation of (\ref{ol_convergence_1}), using (\ref{ol_convergence_2}), that whenever
\begin{equation}
k > \frac{\ln (\frac{\Delta x_{1}^{aug}}{QN})}{\ln \beta} =: N_{1}^{*},
\end{equation}
\noindent $\mathbf{x_{1}^{aug}}(k+N|k) \in \lambda_{1}G_{1} \Rightarrow x_{1}^{aug}(k+N) \in G_{1}$.

Through the same process, one can show that there exists $N_{2}^{*}$ for which
$\tilde{x} \in G_{2}$. Specifically, for the inner loop:
\begin{equation}
\|\mathbf{\tilde{x}}(k+N|k)-\tilde{x}(k+N)\| \leq N (\sum_{j=0}^{i-1} \|A_{2}^{j}B_{2}\|\delta_{u}^{max} + \epsilon_{x_{f}}^{max})\beta^{k},
\end{equation}
\noindent which, after collecting constant terms into one lumped constant $P$ can be rewritten compactly as:
\begin{equation}\label{il_convergence_1}
\|\mathbf{\tilde{x}}(k+N|k)-\tilde{x}(k+N) \leq PN\beta^{k}.
\end{equation}
\noindent Because $\lambda_{2}G_{2} \in G_{2}$, there exists a positive scalar $\Delta \tilde{x}$
such that for any two vectors $\tilde{x}_{a} \in \lambda_{2}G_{2}$ and $\tilde{x}_{b} \in G_{2}$,
$\|\tilde{x}_{a}-\tilde{x}_{b}\|<\Delta \tilde{x}$. To guarantee that $\mathbf{\tilde{x}}(k+N|k) \in \lambda_{1}G_{1} \Rightarrow \tilde{x}(k+N) \in G_{2}$, it suffices to ensure that:
\begin{equation}\label{il_convergence_2}
\|\mathbf{\tilde{x}}(k+N|k)-\tilde{x}(k+N)\| < \Delta \tilde{x}.
\end{equation}

It follows through manipulation of (\ref{il_convergence_1}), using (\ref{il_convergence_2}), that whenever
\begin{equation}
k > \frac{\ln (\frac{\Delta \tilde{x}}{PN})}{\ln \beta} =: N_{2}^{*},
\end{equation}
\noindent $\mathbf{\tilde{x}}(k+N|k) \in \lambda_{2}G_{2} \Rightarrow \tilde{x}(k+N) \in G_{2}$.

Taking $N^{*} \triangleq \max \{N_{1}^{*},N_{2}^{*}\}$ completes the proof.
\begin{flushright} $\Box$
\end{flushright}
\end{pf}
\noindent The proof relies on the fact that the variation in
$\mathbf{v_{des}^{o}}$ and $\mathbf{u^{o}}$ is not only limited, but is also required to decay
over time (through the use of $\beta < 1$ in (\ref{mpc_outer_constraints}) and
(\ref{mpc_inner_constraints})).

\subsection{Overall Stability}
We now state our main result, namely asymptotic stability of the origin
of the overall system, with region of attraction $X$:
\begin{thm} \label{stability} (Asymptotic stability): Under the MPC controller, specified by
(\ref{mpc_outer_law})-(\ref{mpc_inner_cost}), the origin,
$x_{1}^{aug} = 0$, $\tilde{x} = 0$, is asymptotically stable with region of attraction $X$. \end{thm}
\begin{pf} Propositions \ref{terminal_control} and \ref{terminal_control_inexact}
establish the local asymptotic stability of the origin, $x_{1}^{aug}=0, \tilde{x}=0$, under the terminal
control laws, $v_{des}(k)=-K_{1}x_{1}^{aug}(k)$ and $u(k)=u_{t}(v_{des}(k),\tilde{x}(k),x_{f}(k))$.  Because
these terminal control laws are active whenever $x_{1}^{aug} \in G_{1}$ and $\tilde{x} \in G_{2}$, and because
$x_{1}^{aug}(k) \in G_{1}, \tilde{x}(k) \in G_{2} \rightarrow x_{1}^{aug}(k+1) \in G_{1}, \tilde{x}(k+1) \in G_{2}$,
it follows that the origin of the
overall system, $x_{1}^{f} = 0$, $\tilde{x} = 0$, is (locally) asymptotically stable with region of
attraction $\{x_{1}^{aug},\tilde{x}\ : x_{1}^{aug} \in G_{1}, \tilde{x} \in G_{2}\}$.

From Proposition \ref{convergence}, we know that, under the proposed control law, if $\left[\begin{array}{cc} x_{1}^{aug}(0)^{T} & \tilde{x}(0)^{T} \end{array}\right]^{T} \in X$,
then there exists $N^{*}$ for which $x_{1}^{aug}(N^{*}) \in G_{1}$ and
$\tilde{x}(N^{*}) \in G_{2}$.  It follows that $x_{1}^{aug}=0$, $\tilde{x}=0$ has region of attraction $X$.
\begin{flushright} $\Box$
\end{flushright}
\end{pf}
\noindent The proof contains two parts.  First, local asymptotic
stability with region of attraction $\{(x_{1}^{aug},\tilde{x}): x_{1}^{aug} \in G_{1}, \tilde{x} \in G_{2}\}$
is shown by demonstrating that both the inner and outer loop systems are input-to-state stable (ISS)
and the small gain condition is satisfied within this (invariant) region of attraction.
Through the use of MPC, the region of attraction is enlarged to $X$.

\section{Conclusions and Future Work}
In this paper, we reviewed a novel alternative approach to hierarchical MPC that relies
on an inner loop reference model rather than a multi-rate approach for achieving
overall system stability.  This new approach broadens the class of systems for which overall stability of a hierarchical
MPC framework can be guaranteed by allowing the inner closed loop to track the output of a prescribed
reference model rather than requiring the inner loop to reach a steady state at each outer loop step. This paper
presented proofs that were omitted in other works by the authors due to space constraints.

\begin{thebibliography}{xx}  
\bibitem[Vermillion, Menezes, Kolmanovsky (2011)]{hierarchical_vermillion}
C. Vermillion, A. Menezes, I. Kolmanovsky
\newblock Stable Hierarchical Model Predictive Control Using
an Inner Loop Reference Model.
\newblock \emph{Proceedings of the 18th IFAC World Congress}, Milan, Italy, 2011.

\bibitem[Vermillion, Menezes, Kolmanovsky (2013)]{hierarchical_vermillion_journal}
C. Vermillion, A. Menezes, I. Kolmanovsky
\newblock Stable Hierarchical Model Predictive Control Using
an Inner Loop Reference Model and $\lambda$-Contractive Terminal Constraint Sets.
\newblock \emph{Automatica}, Accepted provisionally for publication in December, 2012, Re-submitted in February, 2013.

\bibitem[Falcone et. al.(2008)]{hierarchical_app1}
P. Falcone, F. Borrelli, H. Tseng, J. Asgari, D. Hrovat.
\newblock A Hierarchical Model Predictive Control Framework for Autonomous Ground Vehicles.
\newblock \emph{Proceedings of the American Control Conference}, Seattle, WA, 2008.

\bibitem[Khalil (2001)]{khalil}
H. Khalil
\newblock \emph{Nonlinear Systems}, 3rd Edition
\newblock Prentice Hall, 2001.

\bibitem[Lin, Antsaklis(2004)]{lambda_contractive}
H. Lin and P. Antsaklis.
\newblock A Necessary
and Sufficient Condition for Robust Asymptotic Stabilizability of
Continuous-Time Uncertain Switched Linear Systems.
\newblock \emph{Proceedings
of the IEEE Conference on Decision and Control}, Paradise Island, Bahamas, 2004.

\bibitem[Luenberger(1967)]{miso_ccf}
D. Luenberger.
\newblock Canonical Forms for Linear Multivariable Systems.
\newblock \emph{IEEE Transactions on Automatic Control}, Vol. 12, No. 3, pp. 290-293, 1967.

\bibitem[Luo  et. al.(2004)]{ca2}
Y. Luo, A. Serrani, S. Yurkovich, D. Doman, M. Oppenheimer.
\newblock Model Predictive Dynamic
Control Allocation with Actuator Dynamics.
\newblock \emph{Proceedings of
the American Control Conference}, Boston, MA, 2004.

\bibitem[Luo et. al.(2005)]{ca3}
Y. Luo, A. Serrani, S. Yurkovich, D. Doman, M. Oppenheimer.
\newblock Dynamic Control Allocation with Asymptotic Tracking of Time-Varying Control Input Commands.
\newblock \emph{Proceedings of the American Control Conference}, Portland, OR, 2005.

\bibitem[Luo et. al.(2007)]{ca4}
Y. Luo, A. Serrani, S. Yurkovich, M. Oppenheimer.
\newblock Model Predictive Dynamic Control Allocation Scheme for
Reentry Vehicles.
\newblock \emph{Journal of Guidance, Control, and
Dynamics}, Vol. 30, No. 1, 2007, pp. 100-113.

\bibitem[Mhaskar et. al.(2006)]{application_ex}
P. Mhaskar, N. El-Farra, C. McFall, P. Christofides,
and J. Davis.
\newblock Integrated Fault-Detection and Fault-Tolerant Control of
Process Systems.
\newblock \emph{American Institute of Chemical Engineers Journal}, Vol. 52, pp. 2129-2148, 2006.

\bibitem[Picasso et. al.(2010)]{hierarchical3}
B. Picasso, D. De Vito, R. Scattolini, P. Colaneri.
\newblock An
MPC Approach to the Design of Two-Layer Hierarchical Control Systems.
\newblock \emph{Automatica}, pp. 823-831, 2010.

\bibitem[Scattolini, Colaneri(2007)]{hierarchical1}
R. Scattolini, P. Colaneri.
\newblock Hierarchical
Model Predictive Control.
\newblock H\emph{Proceedings of the IEEE
Conference on Decision and Control}, New Orleans, LA, 2007.

\bibitem[Scattolini et. al.(2008)]{hierarchical2}
R. Scattolini, P, Colaneri, D. Vito.
\newblock A
Switched MPC Approach to Hierarchical Control.
\newblock \emph{Proceedings
of the 17th International Federation of Automatic Control (IFAC) World
Congress}, Seoul, Korea, 2008.

\bibitem[Scattolini(2009)]{hierarchical_review}
R. Scattolini.
\newblock Architectures for Distributed and Hierarchical Model Predictive Control - A Review.
\newblock H\emph{Journal of Process Control}, pp. 723-731, 2009.

\bibitem[Tjonnas, Johansen(2007)]{ca5}
J. Tjonnas, T. Johansen.
\newblock Optimizing Adaptive
Control Allocation with Actuator Dynamics.
\newblock \emph{Proceedings of
the IEEE Conference on Decision and Control}, New Orleans, LA, 2007.

\bibitem[Vermillion et. al.(2007)]{ca6}
C. Vermillion, J. Sun, K. Butts
\newblock Model
Predictive Control Allocation for Overactuated Systems - Stability
and Performance.
\newblock \emph{Proceedings of the IEEE Conference on Decision
and Control}, New Orleans, LA, 2007.

\bibitem[Vermillion et. al.(2009)]{ca7}
C. Vermillion, J. Sun, K. Butts
\newblock Model
Predictive Control Allocation - Design and Experimental Results on
a Thermal Management System.
\newblock \emph{Proceedings of the American Control Conference}, St. Louis, MO, 2009.

\bibitem[Vermillion et. al.(2011)]{ca8}
C. Vermillion, J. Sun, K. Butts
\newblock Predictive Control Allocation for a Thermal Management System
Based on an Inner Loop Reference Model - Design, Analysis, and Experimental
Results.
\newblock \emph{IEEE Transactions on Control Systems Technology}, pp. 772-781, Vol. 19, Issue 4, 2011.
\end{thebibliography}

\end{document}